\documentclass[10pt]{article}
\usepackage{amsmath}
\usepackage{amssymb}
\usepackage{mathrsfs}
\usepackage{amsfonts}
\usepackage{mathrsfs,amscd,amssymb,amsthm,amsmath,bm,graphicx}
\usepackage{graphicx}
\usepackage{cite}
\usepackage{color}
\usepackage{enumitem}
\usepackage{setspace}
\makeatletter

\renewcommand{\@seccntformat}[1]{{\csname the#1\endcsname}{\normalsize.}\hspace{.5em}}
\makeatother

\def \[{\begin{equation}}
\def \]{\end{equation}}

\newtheorem{thm}{Theorem}[section]

\newtheorem{lem}[thm]{Lemma}

\begin{document}
\setlength{\baselineskip}{13pt}
\begin{center}{\Large \bf  Computing the Laplacian spectrum of linear octagonal-quadrilateral networks and its applications
}

\vspace{4mm}

{\large Jia-Bao Liu $^{1,2}$, Zhi-Yu Shi $^{1}$, Ying-Hao Pan$^{1}$, Jinde Cao$^{2,*}$, M. Abdel-Aty$^{3,4}$, Udai Al-Juboori$^{5}$}\vspace{2mm}

{\small $^{1}$School of Mathematics and Physics, Anhui Jianzhu University, Hefei 230601, China\\
$^{2}$School of Mathematics, Southeast University, Nanjing 210096, China\\
$^{3}$Center for Photonics and Smart Materials, Zewail City of Science and Technology, Egypt\\
$^{4}$Mathematics Department, Faculty of Sciences, Sohag University, Egypt\\
$^{5}$College of Arts and Science, Applied Science University, Kingdom of Bahrain}
\vspace{2mm}
\end{center}

\footnotetext{E-mail address: liujiabaoad@163.com, shizhiyuah@163.com, 17318533941@126.com, jdcao@seu.edu.cn,\\
$~~~~~~~~~~~~~~~~~~~~~~~~~~~$ amisaty@gmail.com, udai.al-juboori@asu.edu.bh}

\footnotetext{* Corresponding author.}

 {\noindent{\bf Abstract.}\ \  Let $L_{n}$ denote linear
octagonal-quadrilateral networks. In this paper, we aim to firstly investigate the Laplacian spectrum  on the basis of  Laplacian polynomial of $L_{n}$. Then, by applying the relationship between the coefficients and roots of the
polynomials, the Kirchhoff index and the complexity are determined.

\noindent{\bf Keywords}:  Laplacian matrix; Resistance distance;
Kirchhoff index; Complexity\vspace{2mm}

\noindent{\bf AMS subject classification:} 05C50,\ 05C90}

\section{Introduction}\label{sct1}
\ \ \ \ \ In this article, we only consider simple, undirected and
connected graphs. Suppose $G=(V_{G},E_{G})$ is a graph with vertex set
$V_{G}=\{v_1,v_2,\cdots,v_n\}$ and edge set
$E_{G}=\{e_1,e_2,\cdots,e_m\}$. Let
$D(G)=diag\{d_{1},d_{2},\cdots,d_{n}\}$ be a degree diagonal
matrix, where $d_{i}$ is the degree of $v_{i}$ in $G$. The adjacency matrix $A(G)$ of $G$ is an
$(0,1)$-matrix with order $n$. Then we can get the Laplacian matrix,
which is defined as $L(G)=D(G)-A(G)$. Let
0=$\mu_{1}<\mu_{2}\leq\cdots\leq\mu_{n}$ be the eigenvalues of
$L(G)$. According to the characteristic
polynomial of the matrix $L(G)$, we can get Laplacian spectrum of
$L_{n}$ \cite{C.H}. For more notations and terminologies, one can be referred to \cite{J.}.

At this point, some parameters are introduced. The distance, denoted by $d_{ij}$, is the length of a shortest path between nodes $i$ and $j$, which was named as Wiener index \cite{Wiener,A.D}. This is well-known distance-based topological descriptor, that is
\begin{eqnarray*}
W(G)=\sum_{i<j}d_{ij}.
\end{eqnarray*}

In the electrical network theory, the resistance distance was
firstly proposed by Klein and Randi\'c \cite{D.J}. According to
this concept, we obtain the interpretation of physical
community: the resistance distance between the nodes $i$ and $j$
of the graph $G$ is denoted by $r_{ij}$. One
well-known resistance distance-based parameter called the
Kirchhoff index \cite{D.J.,D.} is given by
$$Kf(G)=\sum_{i<j}r_{ij}.$$

The Kirchhoff index has attracted extensive attentions due to
its wide applications in the fields of physics, chemistry and
others. Despite all that, it is hard to deal with the Kirchhoff index
of complex graphs. Thus, some researchers try to find some new techniques
to compute the Kirchhoff index and obtain its formula. Given an
$n$-vertex graph $G$, Klein and Lov\'asz \cite{Gut,H.Y} proved
independently that
\begin{eqnarray} K f(G)=\sum_{\{u,v\}\subseteq V}
r_{(G)}(u,v)=n\sum_{k=2}^{n}\frac{1}{\mu_{k}} ,
\end{eqnarray}
where $0=\mu_{1}<\mu_{2}\leq\cdots\leq\mu_{n}(n\geq2)$ are the
eigenvalues of $L(G)$.

The number of spanning trees of the graph $G$, also known as the complexity of $G$, is the number of subgraphs that contain all the vertices of $G$ \cite{F.R}. In addition, all those subgraphs must be trees.

According to the decomposition theorem of Laplacian polynomial, Y. Yang et
al., 2008 \cite{Y.} obtained the Laplacian spectrum of linear hexagonal
networks. J. Huang et al. \cite{Huang} got the normalized
Laplacian spectrum of linear hexagonal networks by using the decomposition theorem. Then, the Laplacian spectrum of linear phenylenes
were derived \cite{Y.Peng}. Besides, Z. Zhu and J. Liu
\cite{Zhu} obtained the Laplacian spectrum of generalized
phenylenes. Thus, the extended considerations for calculating the
Laplacian spectrum of linear octagonal-quadrilateral networks are shown in the
following sections.

In the following, we introduce some theorems and notations in Section \ref{sct2}.
Then, we derive the Laplacian spectrum of $L_n$ by using
the relationship between the coefficients and roots in Section \ref{sct3}. An example of the result is given in Section \ref{sct4}. The conclusion is summarized in Section \ref{sct5}.

\section{Preliminary}\label{sct2}
\ \ \ \ \ First, we list some terminology, notations and some
mature consequences in the following.

Given an $n\times n$ matrix, use $M[i_{1},\cdots,i_{k}]$ to be a
submatrix of $M$, which deletes the $i_{1}$-th,
$\cdots$,$i_{k}$-th columns and rows. The characteristic
polynomial of the matrix $M$ is denoted by $P_{M}(x)=det(xI-M)$.

\begin{figure}[htbp]
\centering\includegraphics[width=15cm,height=4cm]{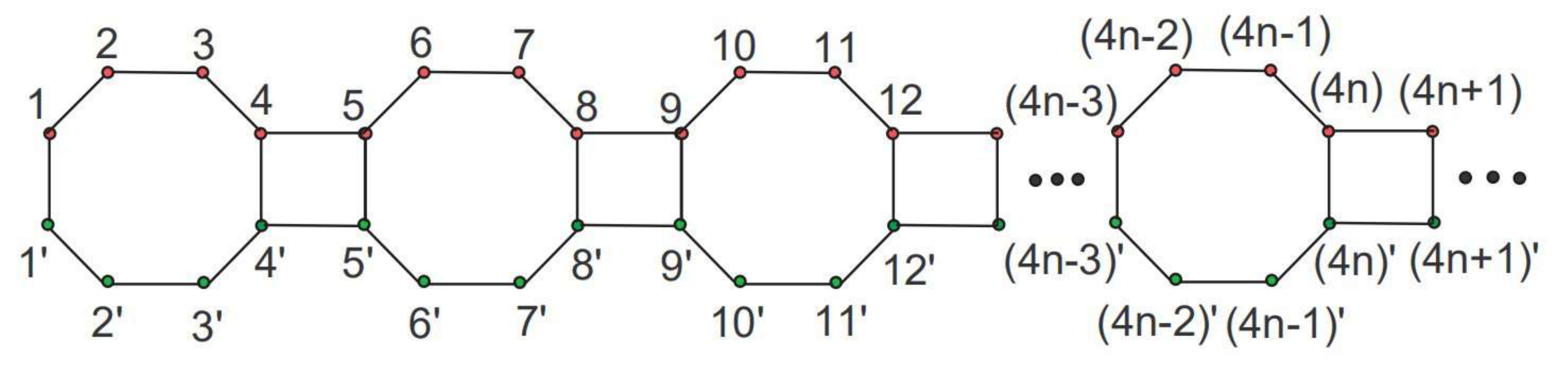}
\caption{Linear octagonal-quadrilateral networks $L_n$.}
\end{figure}

With a suitable labelling of linear octagonal-quadrilateral networks as shown in Figure 1. Evidently,
$|V(L_{n})|=8n+2,|E(L_{n})|=10n+1$. Obviously,
$\pi=(1,1')(2,2')\cdots\big(4n+1,(4n+1)'\big)$ is an automorphism
of $L_{n}$. Set $V_{1}=\{1,2,\cdots,4n+1\},
V_{2}=\{1',2',\cdots,(4n+1)'\}$.

Then $L(L_n)$ can be expressed by the following block matrices.

\begin{equation*}
L(L_{n})=\left(
\begin{array}{cc}
L_{V_{1}V_{1}}& L_{V_{1}V_{2}}\\
L_{V_{2}V_{1}}& L_{V_{2}V_{2}}\\
\end{array}
\right),
\end{equation*}
where
\begin{equation*}
L_{V_{1}V_{1}}=L_{V_{2}V_{2}},~ L_{V_{1}V_{2}}=L_{V_{2}V_{1}}.
\end{equation*}

Let
\begin{equation*}
T=\left(
  \begin{array}{cc}
  \frac{1}{\sqrt{2}}I_{4n+1}& \frac{1}{\sqrt{2}}I_{4n+1}\\
  \frac{1}{\sqrt{2}}I_{4n+1}& -\frac{1}{\sqrt{2}}I_{4n+1}
  \end{array}
\right),
\end{equation*}

then
\begin{equation*}
TL(L_{n})T'=\left(
  \begin{array}{cc}
  L_{A}& 0\\
    0& L_{S}
  \end{array}
\right),
\end{equation*}

where $T'$ is the transposition of $T$ and
\begin{eqnarray}
L_{A}=L_{V_{1}V_{1}}+L_{V_{1}V_{2}},
L_{S}=L_{V_{1}V_{1}}-L_{V_{1}V_{2}}.
\end{eqnarray}

In the block matrix, we can easily get the decomposition theorem
of Laplacian polynomial.

\begin{thm} \textup{\cite{Y. L}}
Assume that $L_{A},L_{S}$ are defined as above. Then we can get
\begin{eqnarray}
P_{L(L_{n})}(x)=P_{L_{A}}(x)P_{L_{S}}(x).
\end{eqnarray}
\end{thm}

\begin{thm} \textup{\cite{W.N}}
If $G$ is the path with $n$ vertices, the eigenvalues of $L(G)$
are $\nu_{i}$, $i=0,1,\cdots,n-1$.
\begin{eqnarray}
\nu_{i}=4sin^{2}(\frac{\pi i}{2n}).
\end{eqnarray}
\end{thm}

\begin{thm} \textup{\cite{F.R}}
If $G$ is a connected graph with $n$ vertices, where $\tau(G)$ is
the complexity of $G$. Then
\begin{eqnarray}
\tau(G)=\frac{1}{n}\prod_{i=2}^{n}\mu_{i}.
\end{eqnarray}
\end{thm}

\begin{thm} \textup{\cite{K.H.}}
(Matrix Tree Theorem) If $G$ is a connected graph with $n$ vertices and $L$ be the Laplacian matrix of $G$. Then the complexity of $G$ is
\begin{eqnarray}
\tau(G)=det(L(i)),
\end{eqnarray}
where $i=1,2,\cdots,n$.
\end{thm}

In the following, a flowchart is given according to the steps we have processed, which helps to understand the proposed approach. The explanations of these notations that appear in the flowchart are describes in Section \ref{sct3}.

\begin{figure}[htbp]
\centering\includegraphics[width=10cm,height=5cm]{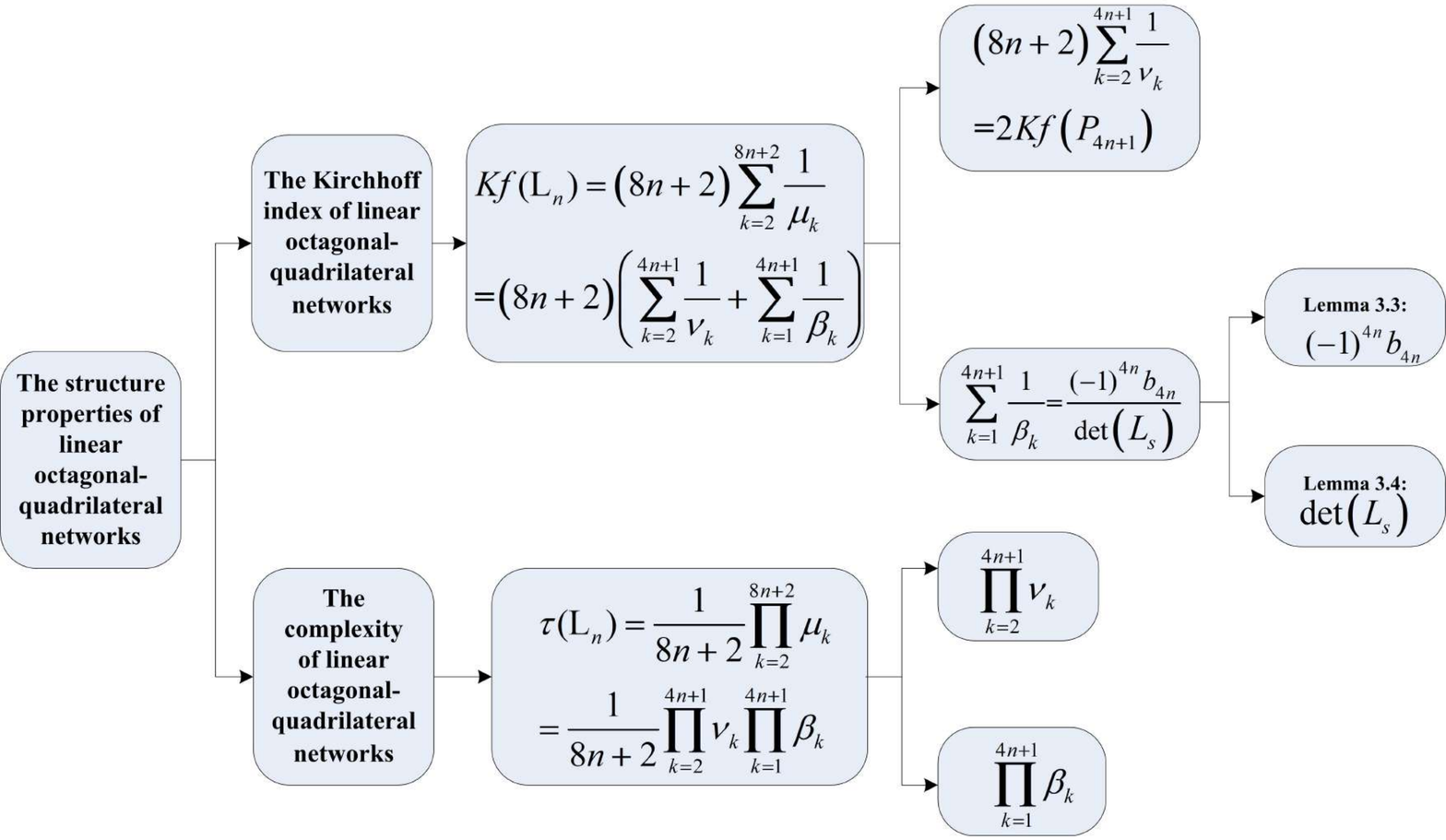}
\caption{The steps processed in this paper.}
\end{figure}

\section{Kirchhoff index and the complexity of $L_{n}$}\label{sct3}
\ \ \ \ \ In this section, we get the general formula of
the Kirchhoff index and the complexity. Using the decomposition
theorem, we can get

\begin{eqnarray*}
 L_{V_1 V_1}&=&
\left(
  \begin{array}{cccccccccccccc}
    2 & -1 & & & & & & & & & & & &\\
    -1 & 2 & -1 & & & & & & & & & & &\\
    & -1 & 2 & -1 & & & & & & & & & &\\
    & & -1 & 3 & -1 & & & & & & & & &\\
    & & & -1 & 3 & -1 & & & & & & & &\\
    & & & & -1 & 2 & -1 & & & & & & &\\
    & & & & & -1 & 2 & -1 & & & & & &\\
    & & & & & & -1 & 3 & -1 & & & & &\\
    & & & & & & & & \ddots & & & & &\\
    & & & & & & & & -1 & 3 & -1 & & &\\
    & & & & & & & & & -1 & 2 & -1 & &\\
    & & & & & & & & & & -1 & 2 & -1 &\\
    & & & & & & & & & & & -1 & 3 & -1\\
    & & & & & & & & & & & & -1 & 2\\
  \end{array}
\right)_{(4n+1)\times (4n+1)},
\end{eqnarray*}

\begin{eqnarray*}
L_{V_1
V_2}=diag\{-1,0,0-1,-1,0,0,-1,\cdots,-1,0,0,-1,-1\}_{(4n+1)}.
\end{eqnarray*}

Through Theorem 2.1, one gets the Laplacian spectrum, which are composed of the eigenvalues of $L_{A}$ and $L_{S}$ of $L_{n}$ as follows.

\begin{eqnarray*}
 L_A&=&
\left(
  \begin{array}{cccccccccc}
    1 & -1 & & & & & & & &\\
    -1 & 2 & -1 & & & & & & &\\
    & -1 & 2 & -1 & & & & & &\\
    & & -1 & 2 & -1 & & & & &\\
    & & & -1 & 2 & -1 & & & &\\
    & & & & & \ddots & & & &\\
    & & & & & -1 & 2 & -1 & &\\
    & & & & & & -1 & 2 & -1 &\\
    & & & & & & & -1 & 2 & -1\\
    & & & & & & & & -1 & 1\\
  \end{array}
\right)_{(4n+1)\times (4n+1)},
\end{eqnarray*}

\begin{eqnarray*}
 L_S&=&
\left(
  \begin{array}{cccccccccccccc}
    3 & -1 & & & & & & & & & & & &\\
    -1 & 2 & -1 & & & & & & & & & & &\\
    & -1 & 2 & -1 & & & & & & & & & &\\
    & & -1 & 4 & -1 & & & & & & & & &\\
    & & & -1 & 4 & -1 & & & & & & & &\\
    & & & & -1 & 2 & -1 & & & & & & &\\
    & & & & & -1 & 2 & -1 & & & & & &\\
    & & & & & & -1 & 4 & -1 & & & & &\\
    & & & & & & & & \ddots & & & & &\\
    & & & & & & & & -1 & 4 & -1 & & &\\
    & & & & & & & & & -1 & 2 & -1 & &\\
    & & & & & & & & & & -1 & 2 & -1 &\\
    & & & & & & & & & & & -1 & 4 & -1\\
    & & & & & & & & & & & & -1 & 3\\
  \end{array}
\right)_{(4n+1)\times (4n+1)}.
\end{eqnarray*}

By calculation, we find that $L_{A}$ is obviously the Laplacian
matrix of the path $P_{4n+1}$.

Assume that
$0=\nu_{1}<\nu_{2}\leq\nu_{3}\leq\cdots\leq\nu_{4n+1}$ are the
roots of polynomial $P_{L_{A}}(x)=0$, and
$0<\beta_{1}\leq\beta_{2}\leq\beta_{3}\leq\cdots\leq\beta_{4n+1}$ are the
roots of polynomial $P_{L_{S}}(x)=0$. Therefore, we have
\begin{eqnarray}
K
f(L_{n})=(8n+2)\Big(\sum_{k=2}^{4n+1}\frac{1}{\nu_{k}}+\sum_{k=1}^{4n+1}\frac{1}{\beta_{k}}\Big)=2K
f(P_{4n+1})+(8n+2)\sum_{k=1}^{4n+1}\frac{1}{\beta_{k}}.
\end{eqnarray}

By applying the famous formula $Kf(P_{n})=\frac{n^{3}-n}{6}$, we
can get
\begin{eqnarray}
2K f(P_{4n+1})=\frac{8n(4n+1)(2n+1)}{3}.
\end{eqnarray}

Thus, we only need to calculate the latter part
$\sum_{k=1}^{4n+1}\frac{1}{\beta_{k}}$.

Let
$$P_{L_{S}}(x)=det(xI-L_{S})=x^{4n+1}+b_{1}x^{4n}+\cdots+b_{4n}x+b_{4n+1}.$$

By the relationship between coefficients and roots of the
characteristic polynomial $P_{L_{S}}(x)$, we can get
\begin{eqnarray}
\sum_{k=1}^{4n+1}\frac{1}{\beta_{k}}=\frac{(-1)^{4n}b_{4n}}{(-1)^{4n+1}b_{4n+1}}=\frac{(-1)^{4n}b_{4n}}{det(L_{S})}.
\end{eqnarray}

Obviously, we obtain the
$P_{L_{S}}(0)=(-1)^{4n+1}det(L_{S})=b_{4n+1}$ and the
$(-1)^{4n}b_{4n}$ to be equal to sum of the principal minors of all
$4n$ columns and $4n$ rows of $L_{S}$. Therefore, we make the
$k$-th order principal submatrix to be $C_{k}$, which is composed of
the first $k$ columns and rows of $L_{S}$, $k=1,2,\cdots,4n$. Let
$c_{k}=det(C_{k})$. Then
$c_{1}=3,c_{2}=5,c_{3}=7,c_{4}=23,c_{5}=85,c_{6}=147,c_{7}=209,c_{8}=689$
and for $i\geq 1$,
\begin{eqnarray*}
\begin{cases}
c_{4i+1}=4c_{4i}-c_{4i-1};\\
c_{4i+2}=2c_{4i+1}-c_{4i};\\
c_{4i+3}=2c_{4i+2}-c_{4i+1};\\
c_{4i+4}=4c_{4i+3}-c_{4i+2}.
\end{cases}
\end{eqnarray*}

In the following, one gets the further terms of the above recurrence formulae.
\begin{eqnarray*}
\begin{cases}
c_{4i+1}=30c_{4i-3}-c_{4i-7};\\
c_{4i+2}=30c_{4i-2}-c_{4i-6};\\
c_{4i+3}=30c_{4i-1}-c_{4i-5};\\
c_{4i+4}=30c_{4i}-c_{4i-4}.
\end{cases}
\end{eqnarray*}

For convenience, let
$A$ and $B$ be $15+4\sqrt14$ and $15-4\sqrt14$. Throughout the
rest of the context, we omit to introduce the definitions of $A$
and $B$, if there are no confusions.

\begin{lem}
Let ${c_{4i+1}}_{(i\geq0)}$, ${c_{4i+2}}_{(i\geq0)}$,
${c_{4i+3}}_{(i\geq0)}$ and ${c_{4i+4}}_{(i\geq0)}$ be the sequences defined as above, for $0\leq
i\leq n-1$,
\begin{eqnarray*}
\begin{cases}
c_{4i+1}=\frac{(10+3\sqrt{14})A^{i}-(10-3\sqrt{14})B^{i}}{2\sqrt{14}};\\\\
c_{4i+2}=\frac{(18+5\sqrt{14})A^{i}-(18-5\sqrt{14})B^{i}}{2\sqrt{14}};\\\\
c_{4i+3}=\frac{(26+7\sqrt{14})A^{i}-(26-7\sqrt{14})B^{i}}{2\sqrt{14}};\\\\
c_{4i+4}=\frac{(86+23\sqrt{14})A^{i}-(86-23\sqrt{14})B^{i}}{2\sqrt{14}}.
\end{cases}
\end{eqnarray*}
\end{lem}

Now we make the $k$-th order principal submatrix to be $E_{k}$, which
is composed of the last $k$ columns and rows of $L_{S}$,
$k=1,2,\cdots,4n$. Let $e_{k}=det(E_{k})$. Then
$e_{1}=3,e_{2}=11,e_{3}=19,e_{4}=27,e_{5}=89,e_{6}=329,e_{7}=569,e_{8}=809$
and for $i\geq1$,
\begin{eqnarray*}
\begin{cases}
e_{4i+1}=4e_{4i}-e_{4i-1};\\
e_{4i+2}=4e_{4i+1}-e_{4i};\\
e_{4i+3}=2e_{4i+2}-e_{4i+1};\\
e_{4i+4}=2e_{4i+3}-e_{4i+2}.
\end{cases}
\end{eqnarray*}

Then, we can get
\begin{eqnarray*}
\begin{cases}
e_{4i+1}=30e_{4i-3}-e_{4i-7};\\
e_{4i+2}=30e_{4i-2}-e_{4i-6};\\
e_{4i+3}=30e_{4i-1}-e_{4i-5};\\
e_{4i+4}=30e_{4i}-e_{4i-4}.
\end{cases}
\end{eqnarray*}

Using the above method, we can solve the infinite sequences
${e_{4i}}_{(i\geq1)}$
(resp.${e_{4i+1}}_{(i\geq0)}$, ${e_{4i+2}}_{(i\geq0)}$, ${e_{4i+3}}_{(i\geq0)}$)
and get the general terms.

\begin{lem}
Let ${e_{4i+1}}_{(i\geq0)}$, ${e_{4i+2}}_{(i\geq0)}$,
${e_{4i+3}}_{(i\geq0)}$ and ${e_{4i+4}}_{(i\geq0)}$ be the sequences defined as above, for $0\leq
i\leq n-1$,
\begin{eqnarray*}
\begin{cases}
e_{4i+1}=\frac{(11+3\sqrt{14})A^{i}-(11-3\sqrt{14})B^{i}}{2\sqrt{14}};\\\\
e_{4i+2}=\frac{(41+11\sqrt{14})A^{i}-(41-11\sqrt{14})B^{i}}{2\sqrt{14}};\\\\
e_{4i+3}=\frac{(71+19\sqrt{14})A^{i}-(71-19\sqrt{14})B^{i}}{2\sqrt{14}};\\\\
e_{4i+4}=\frac{(101+27\sqrt{14})A^{i}-(101-27\sqrt{14})B^{i}}{2\sqrt{14}}.
\end{cases}
\end{eqnarray*}
\end{lem}

\begin{lem} $(-1)^{4n}b_{4n}=\frac{1}{196}[(98+33\sqrt{14})A^{n}+(98-33\sqrt{14})B^{n}]
+\frac{9n}{14}[(4+\sqrt{14})A^{n}+(4-\sqrt{14})B^{n}]$.
\end{lem}

\noindent{\bf Proof.} Since numeric expression $(-1)^{4n}b_{4n}$ is equal to sum of the
principal minors of all $4n$ columns and $4n$ rows of $L_{S}$, we
can get

\begin{eqnarray*}
(-1)^{4n}b_{4n}&=&\sum_{i=1}^{4n+1}det(L_{S}[i])\\
&=&\sum_{k=0}^{n}detL_{S}[4k+1]+\sum_{k=0}^{n-1}detL_{S}[4k+2]+\sum_{k=0}^{n-1}detL_{S}[4k+3]+\sum_{k=0}^{n-1}detL_{S}[4k+4]\\
&=&c_{4n}+e_{4n}+\sum_{k=1}^{n-1}c_{4k}e_{4n-4k}+\sum_{k=0}^{n-1}c_{4k+1}e_{4n-4k-1}+\sum_{k=0}^{n-1}c_{4k+2}e_{4n-4k-2}+\sum_{k=0}^{n-1}c_{4k+3}e_{4n-4k-3}\\
&=&c_{4n}+e_{4n}+\sum_{k=1}^{n-1}c_{4k}e_{4(n-k)}+\sum_{k=0}^{n-1}c_{4k+1}e_{4(n-k-1)+3}+\sum_{k=0}^{n-1}c_{4k+2}e_{4(n-k-1)+2}\\
&&+\sum_{k=0}^{n-1}c_{4k+3}e_{4(n-k-1)+1}.
\end{eqnarray*}

For convenience, we set
\begin{eqnarray*}
\begin{cases}
D_{1}=\sum\limits_{k=1}^{n-1}c_{4k}e_{4(n-k)};\\\\
D_{2}=\sum\limits_{k=0}^{n-1}c_{4k+1}e_{4(n-k-1)+3};\\\\
D_{3}=\sum\limits_{k=0}^{n-1}c_{4k+2}e_{4(n-k-1)+2};\\\\
D_{4}=\sum\limits_{k=0}^{n-1}c_{4k+3}e_{4(n-k-1)+1}.
\end{cases}
\end{eqnarray*}

Then

\begin{eqnarray*}
D_{1}&=&\sum_{k=1}^{n-1}c_{4k}e_{4(n-k)}\\
&=&\sum_{k=1}^{n-1}\Big[\frac{(86+23\sqrt{14})A^{k-1}-(86-23\sqrt{14})B^{k-1}}{2\sqrt{14}}\cdot\frac{(101+27\sqrt{14})A^{n-k-1}-(101-27\sqrt{14})B^{n-k-1}}{2\sqrt{14}}\Big]\\
&=&\frac{n-1}{56}[(17380+4645\sqrt{14})A^{n-2}+(17380-4645\sqrt{14})B^{n-2}]+\frac{\sqrt{14}}{392}(A^{n-1}-B^{n-1}).
\end{eqnarray*}

\begin{eqnarray*}
D_{2}&=&\sum_{k=0}^{n-1}c_{4k+1}e_{4(n-k-1)+3}\\
&=&\sum_{k=0}^{n-1}\Big[\frac{(10+3\sqrt{14})A^{k}-(10-3\sqrt{14})B^{k}}{2\sqrt{14}}\cdot\frac{(71+19\sqrt{14})A^{n-k-1}-(71-19\sqrt{14})B^{n-k-1}}{2\sqrt{14}}\Big]\\
&=&\frac{n}{56}[(1508+403\sqrt{14})A^{n-1}+(1508-403\sqrt{14})B^{n-1}]+\frac{11\sqrt{14}}{392}(A^{n}-B^{n}).
\end{eqnarray*}

\begin{eqnarray*}
D_{3}&=&\sum_{k=0}^{n-1}c_{4k+2}e_{4(n-k-1)+2}\\
&=&\sum_{k=0}^{n-1}\Big[\frac{(18+5\sqrt{14})A^{k}-(18-5\sqrt{14})B^{k}}{2\sqrt{14}}\cdot\frac{(41+11\sqrt{14})A^{n-k-1}-(41-11\sqrt{14})B^{n-k-1}}{2\sqrt{14}}\Big]\\
&=&\frac{n}{56}[(1508+403\sqrt{14})A^{n-1}+(1508-403\sqrt{14})B^{n-1}]+\frac{4\sqrt{14}}{392}(A^{n}-B^{n}).
\end{eqnarray*}

\begin{eqnarray*}
D_{4}&=&\sum_{k=0}^{n-1}c_{4k+3}e_{4(n-k-1)+1}\\
&=&\sum_{k=0}^{n-1}\Big[\frac{(26+7\sqrt{14})A^{k}-(26-7\sqrt{14})B^{k}}{2\sqrt{14}}\cdot\frac{(11+3\sqrt{14})A^{n-k-1}-(11-3\sqrt{14})B^{n-k-1}}{2\sqrt{14}}\Big]\\
&=&\frac{n}{56}[(580+155\sqrt{14})A^{n-1}+(580-155\sqrt{14})B^{n-1}]+\frac{\sqrt{14}}{392}(A^{n}-B^{n}).
\end{eqnarray*}

Therefore, we can have
$(-1)^{4n}b_{4n}=c_{4n}+e_{4n}+D_{1}+D_{2}+D_{3}+D_{4}$. Then
\begin{eqnarray*}
(-1)^{4n}b_{4n}&=&\frac{1}{196}[(98+33\sqrt{14})A^{n}+(98-33\sqrt{14})B^{n}]+\frac{9n}{14}[(4+\sqrt{14})A^{n}+(4-\sqrt{14})B^{n}].
\end{eqnarray*}

This completes the proof. \hfill\rule{1ex}{1ex}

\begin{lem}
$det(L_{S})=\frac{(116+31\sqrt{14})A^{n-1}-(116-31\sqrt{14})B^{n-1}}{\sqrt{14}}.$
\end{lem}

\noindent{\bf Proof.} By expanding $L_{s}$ with regard to the last row, we have
\begin{eqnarray*}
det(L_{S})&=&3c_{4n}-c_{4(n-1)+3}\\
&=&\frac{3[(86+23\sqrt{14})A^{n-1}-(86-23\sqrt{14})B^{n-1}]}{2\sqrt{14}}-\frac{[(26+7\sqrt{14})A^{n-1}-(26-7\sqrt{14})B^{n-1}]}{2\sqrt{14}}\\
&=&\frac{(116+31\sqrt{14})A^{n-1}-(116-31\sqrt{14})B^{n-1}}{\sqrt{14}}.
\end{eqnarray*}

This completes the proof. \hfill\rule{1ex}{1ex}

Together with formulas (3.7)~-~(3.9) and lemmas 3.3~-~3.4, one can
get the following theorem.

\begin{thm}
For linear octagonal-quadrilateral networks $L_{n}$,
$$Kf(L_{n})=\frac{8n(4n+1)(2n+1)}{3}+\frac{(8n+2)(-1)^{4n}b_{4n}}{detL_{S}},$$
where
\begin{eqnarray*}
(-1)^{4n}b_{4n}&=&\frac{1}{196}[(98+33\sqrt{14})A^{n}+(98-33\sqrt{14})B^{n}]+\frac{9n}{14}[(4+\sqrt{14})A^{n}+(4-\sqrt{14})B^{n}],\\
detL_{S}
&=&\frac{(116+31\sqrt{14})A^{n-1}-(116-31\sqrt{14})B^{n-1}}{\sqrt{14}}.
\end{eqnarray*}
\end{thm}

The Kirchhoff indices of $L_{n}$ are shown in Table 1, in which $n$ is from 1
to 15.

\begin{table}[htbp]
\setlength{\abovecaptionskip}{0.05cm} \centering\vspace{.3cm}
\caption{Kirchhoff indices from $L_{1}$ to $L_{15}$}
\begin{tabular}{c|c|c|c|c|c|c|c|c|c}
  \hline
  $G$ & $Kf(G)$ & $G$ & $Kf(G)$ & $G$ & $Kf(G)$ & $G$ & $Kf(G)$ & $G$ & $Kf(G)$ \\
  \hline
  $L_{1}$ & $69.52$ & $L_{4}$ & $1977.70$ & $L_{7}$ & $9128.25$ & $L_{10}$ & $24977.17$ & $L_{13}$ & $52980.46$ \\
  $L_{2}$ & $336.42$ & $L_{5}$ & $3608.06$ & $L_{8}$ & $13274.07$ & $L_{11}$ & $32790.45$ & $L_{14}$ & $65613.20$
  \\
  $L_{3}$ & $929.82$ & $L_{6}$ & $5948.91$ & $L_{9}$ & $18514.38$ & $L_{12}$ & $42082.21$ & $L_{15}$ & $ 80108.42$
  \\
  \hline
\end{tabular}
\end{table}

The explicit formula for the complexity of $L_{n}$ is in the following.

\begin{thm}
For octagonal-quadrilateral networks $L_{n}$,
$$\tau(L_{n})=\frac{4^{4n}[(116+31\sqrt{14})A^{n-1}-(116-31\sqrt{14})B^{n-1}]}{\sqrt{14}(8n+2)}
\prod_{i=2}^{4n+1}sin^{2}\Big(\frac{\pi(k-1)}{8n+2}\Big).$$
\end{thm}

\noindent{\bf Proof.} According to Theorem 2.2, we get the
eigenvalues of $L_{A}$, which are
$\nu_{k}=4sin^{2}(\frac{\pi(k-1)}{8n+2})$ ($k=1,2,\cdots,4n+1$).

Together with Lemma 3.4 and Theorem 2.3, we know that
\begin{eqnarray*}
\tau(L_{n})&=&\frac{1}{8n+2}\prod_{k=2}^{4n+1}\nu_{k}\cdot\prod_{k=1}^{4n+1}\beta_{k}\\
&=&\frac{det(L_{s})}{8n+2}\prod_{k=2}^{4n+1}4sin^{2}\Big(\frac{\pi(k-1)}{8n+2}\Big)\\
&=&\frac{4^{4n}[(116+31\sqrt{14})A^{n-1}-(116-31\sqrt{14})B^{n-1}]}{\sqrt{14}(8n+2)}
\prod_{k=2}^{4n+1}sin^{2}\Big(\frac{\pi(k-1)}{8n+2}\Big).
\end{eqnarray*}

This completes the proof. \hfill\rule{1ex}{1ex}

The complexity of $L_{n}$ is shown Table 2, in which $n$ is from 1 to 12.

\begin{table}[htbp]
\setlength{\abovecaptionskip}{0.05cm}
 \centering \vspace{.3cm}
\caption{The complexity from $L_{1}$ to $L_{12}$}
\begin{tabular}{c|c|c|c|c|c|c|c}
  \hline
  $G$ & $\tau(G)$ & $G$ & $\tau(G)$ & $G$ & $\tau(G)$ & $G$ & $\tau(G)$ \\
  \hline
  $L_{1}$ & $31$ & $L_{4}$ & $834241$ & $L_{7}$ & $22449425279$ & $L_{10}$ & $604114033423649$  \\
  $L_{2}$ & $929$ & $L_{5}$ & $24999391$ & $L_{8}$ & $672733610881$ & $L_{11}$ & $18103261443808319$
  \\
  $L_{3}$ & $27839$ & $L_{6}$ & $749147489$ & $L_{9}$ & $20159558901151$ & $L_{12}$ & $542493729280825921$
  \\
  \hline
\end{tabular}
\end{table}

\section{Example}\label{sct4}
\ \ \ \ \ For example, we calculate the Kirchhoff index and the complexity for $L_{1}$. First, we use the decomposition theorem to obtain the block matrix. Then, the relationship between the coefficients and roots derives the explicit formulas for the Kirchhoff index and the complexity. Thus, we can obtain
$L_{V_1 V_1}$ and $L_{V_2 V_2}$, as follows.

\begin{eqnarray*}
 L_{V_1 V_1}&=&
\left(
  \begin{array}{ccccc}
    2 & -1 & & & \\
    -1 & 2 & -1 & & \\
    & -1 & 2 & -1 & \\
    & & -1 & 3 & -1 \\
    & & & -1 & 2 \\
 \end{array}
 \right),
\end{eqnarray*}

\begin{eqnarray*}
 L_{V_2 V_2}&=&
\left(
  \begin{array}{ccccc}
    -1 & & & &\\
    & 0 & & &\\
    & & 0 & &\\
    & & & -1 &\\
    & & & & -1\\
   \end{array}
\right).
\end{eqnarray*}

According to the decomposition theorem, one gets the two special
matrices, $L_{A}$ and $L_{S}$ with order $5$, as
follows

\begin{eqnarray*}
 L_A&=&
\left(
  \begin{array}{ccccc}
    1 & -1 & & &\\
    -1 & 2 & -1 & &\\
    & -1 & 2 & -1 &\\
    & & -1 & 2 & -1\\
    & & & -1 & 1\\
   \end{array}
\right),
\end{eqnarray*}

\begin{eqnarray*}
 L_S&=&
\left(
  \begin{array}{ccccc}
    3 & -1 & & &\\
    -1 & 2 & -1 & &\\
    & -1 & 2 & -1 &\\
    & & -1 & 4 & -1\\
    & & & -1 & 3\\
  \end{array}
\right).
\end{eqnarray*}

Obviously, $L_{A}$ is the Laplacian matrix of path $P_{5}$. Let
$0=\nu_{1}<\nu_{2}\leq\nu_{3}\leq\nu_{4}\leq\nu_{5}$ be the
roots of characteristic polynomial $P_{L_{A}}(x)=0$, and
$0<\beta_{1}\leq\beta_{2}\leq\beta_{3}\leq\beta_{4}\leq\beta_{5}$ be the
roots of characteristic polynomial $P_{L_{S}}(x)=0$.

Thus, we can obtain
\begin{eqnarray*}
 Kf(G)=2Kf(P_{5})+10\sum_{k=1}^{5}\frac{1}{\beta_{k}}=\frac{2155}{31},
\end{eqnarray*}
\begin{eqnarray*}
\tau(G)=\frac{1}{10}\prod_{k=2}^{5}\nu_{k}\prod_{k=1}^{5}\beta_{k}=31.
\end{eqnarray*}

\noindent{\bf Remark.} According to Theorem 2.4, we can calculate the complexity of $L_{1}$. Meanwhile, it is found that the complexity is equal to that result which has been calculated by Laplacian spectrum.

\section{Conclusion}\label{sct5}
\ \ \ \ \ We mainly use the decomposition theorem of Laplacian polynomial. The relationship between coefficients and roots of
the polynomial is a necessary method for us to arrive the
Kirchhoff index and the complexity of octagonal-quadrilateral networks.
In addition, this method can be applied to other graphs.

\section*{Author contribution}
\ \ \ \ \ Funding acquisition, J.-B Liu; Methodology, J.-B Liu, Jinde Cao; Formal analysis, J.-B Liu, Zhi-Yu Shi; Data curation, Writing-original draft, Zhi-Yu Shi, Ying-Hao Pan, M. Abdel-Aty, Udai Al-Juboori.

\section*{Funding}

\ \ \ \ \ The work was partly supported by China Postdoctoral
Science Foundation (No. 2017M621579), Postdoctoral Science
Foundation of Jiangsu Province (No. 1701081B)  and Project of
Anhui Jianzhu University (No. 2016QD116, 2017dc03 and 2017QD20).



\end{document}